\documentclass[11 pt]{amsart}     
\usepackage{amsthm}
\usepackage{amssymb}
\usepackage{latexsym}
\usepackage{amsfonts}
\usepackage{amsmath}
\usepackage{amstext}
\usepackage{amsrefs}
\newtheorem{theorem}{Theorem}

\newtheorem{proposition}[theorem]{Proposition}
\newcommand{\RR}{\mathbb{R}}

\newcommand{\E}{\mathrm{e}}

\begin{document}

\title[Isoperimetric comparison for Ricci flow on the two-sphere]{Curvature bounds by isoperimetric comparison for normalized Ricci flow on the two-sphere}
\thanks{The research of both authors was supported by Discovery Grant DP0556211 of the Australian Research Council.  The second author's research was assisted by an Australian Postgraduate Award.}

\author{Ben Andrews}


\address{MSI, ANU,
              ACT 0200 Australia}
              
              \email{Ben.Andrews@anu.edu.au}           
\author{Paul Bryan}
\address{MSI, ANU, ACT 0200 Australia}
\email{Paul.Bryan@anu.edu.au}

\date{August 25, 2009}

\begin{abstract}
We prove a comparison theorem for the isoperimetric profiles of solutions of the normalized Ricci flow on the two-sphere:  If the isoperimetric profile of the initial metric is greater than that of some positively curved axisymmetric metric, then the inequality remains true for the isoperimetric profiles of the evolved metrics.    We apply this using the Rosenau solution as the model metric to deduce sharp time-dependent curvature bounds for arbitrary solutions of the normalized Ricci flow on the two-sphere.  This gives a simple and direct proof of convergence to a constant curvature metric without use of any blowup or compactness arguments, Harnack estimates, or any classification of behaviour near singularities.
\end{abstract}

\maketitle

\section{Introduction}

The Ricci flow is a nonlinear parabolic evolution for Riemannian metrics which has been applied widely and with great effect in recent years, most notably in the work of Perel'man \cites{Pe1,Pe2,Pe3} proving the Poincar\'e and geometrization conjectures, and in work of B\"ohm and Wilking \cite{BW} and Brendle and Schoen \cites{BS1,BS2,Bre} in proving a series of long-conjectured sphere theorems.  The Ricci flow was introduced by Richard Hamilton who applied the flow to Riemannian metrics of positive Ricci curvature on three-dimensional manifolds \cite{Ham82}.  

The Ricci flow on surfaces is of a rather different character than in higher dimensions, partly because there is only a single curvature function at each point, so that methods involving pointwise comparisons of different parts of the curvature tensor give no information here.  In particular the tensor maximum principle arguments which feature strongly in higher dimensions are of much less utility in two dimensions, and very different proofs are required.  

The Ricci flow on surfaces of genus at least one can be dealt with without great difficulty 
\cites{Ham2DRF1} by introducing a potential function which allows direct estimates on the curvature using the maximum principle.  However the genus zero case is much more subtle, and was resolved only in several stages:  Hamilton \cite{Ham2DRF1} proved that the flow beginning with a positively curved metric will converge to constant curvature, and Bennett Chow \cite{ChowRFS2} proved the result for arbitrary metrics on with genus zero, by making use of a differential Harnack inequality and an entropy estimate to show that the curvature eventually becomes positive.

The method we employ in this paper is inspired by more recent work of Richard Hamilton \cite{Ham2DRF2}, in which isoperimetric bounds were derived to rule out the possibility of a `type II' singularity.  This, together with his earlier work in classifying both Type I and Type II singularities, led a to a proof of convergence to constant curvature for arbitrary metrics on the two-sphere.  In this paper we show that a modification to his argument directly implies bounds on curvature for solutions of the normalized Ricci flow on the two-sphere, sufficiently strong that they give exponential convergence to a constant curvature metric without requiring any further tools:  No Harnack estimates or entropy estimates, compactness theorems, or classifications of type I or type II solutions are required.

Our result amounts to a comparison theorem, allowing the isoperimetric profile of an arbitrary solution of normalized Ricci flow on the two-sphere to be compared with that of an arbitrary axially symmetric positively curved solution.  Our explicit bounds are obtained by applying this comparison with the explicit solution known as the Rosenau solution, and consequently the bounds we obtain are sharp, and attained exactly on the Rosenau solution.

\section{Notation and preliminary results}

We denote by $K(x)$ the Gaussian curvature of a surface $M$ at $x$, given by $R(e_1,e_2,e_1,e_2)$ for any orthonormal basis for the tangent space at $x$ and hence equal to half of the scalar curvature.  The normalized Ricci flow is the equation
\begin{equation}\label{eq:nRF}
\frac{\partial g}{\partial t} = -2(K-k)g,
\end{equation}
where $k$ is the average of the scalar curvature: $k=\frac{1}{|M|}\int_M K d\mu(g)$, which by the Gauss-Bonnet theorem is equal to $\frac{4\pi}{|M|}$ in the case of a metric on the 2-sphere.  Observing that the normalized Ricci flow preserves the area of the surface, we scale the initial metric so that $|M|=4\pi$, and the Ricci flow equation becomes $\frac{\partial g}{\partial t}=-2(K-1)g$.

We make use of an elementary existence result proved as in \cite{Ham82}*{Theorem 14.1}:  For any smooth normalized initial metric there exists a unique solution of equation \eqref{eq:nRF}, which continues to exist and remains smooth as long as $K$ remains bounded.

For a compact Riemannian surface $M$, the \emph{isoperimetric profile} $h: (0,1)\to\RR_+$ of $M$ is defined by
$$
h(\xi)=\inf\left\{|\partial\Omega|:\ \Omega\subset M,\ |\Omega|=\xi|M|\right\}.
$$
For each $\xi\in(0,1)$ the isoperimetric ratio is attained by a domain $\Omega\subset M$ with boundary given by a smooth embedded curve with constant geodesic curvature.  We will need the following observation about the topology of the optimal curves:

\begin{theorem}\label{thm:curvetopology}
Let $(M,g)$ be a compact Riemannian surface with isoperimetric profile $h$.  Let $\varphi$ be a strictly concave positive function on $(0,1)$.  If $h(\xi)\geq \varphi(\xi)$ for every $\xi$, and $h(\xi_0)=\varphi(\xi_0)$ for some $\xi_0\in(0,1)$, then the optimal domain $\Omega$ with $|\Omega|=\xi_0|M|$ which satisfies $|\partial\Omega|=\varphi(\xi_0)$ is connected and has connected complement.
In particular, if $M=S^2$ then $\Omega$ is simply connected.
\end{theorem}

\begin{proof}
We prove that $\Omega$ and its complement are connected:  If not, then we can write $\Omega=\Omega_1\cup\Omega_2$, and we have $|\Omega|=|\Omega_1|+|\Omega_2|$ and $|\partial\Omega|=|\partial\Omega_1|+|\partial\Omega_2|$, and all of these are non-zero.  But then we have
\begin{align*}
\varphi\left(\frac{|\Omega_1|}{|M|}\right)+\varphi\left(\frac{|\Omega_2|}{|M|}\right)&\leq |\partial\Omega_1|+|\partial\Omega_2|\\
&=|\partial\Omega|\\
&=\varphi\left(\frac{|\Omega|}{|M|}\right)\\
&=\varphi\left(\frac{|\Omega_1|}{|M|}+\frac{|\Omega_2|}{|M|}\right)\\
&<\varphi\left(\frac{|\Omega_1|}{|M|}\right)+\varphi\left(\frac{|\Omega_2|}{|M|}\right),
\end{align*}
where we used strict concavity and non-negativity of $\varphi$ in the last inequality (observing that the arguments of the two terms in the last expression are strictly between $0$ and $1$).  This is a contradiction, so $\Omega$ is connected.  Since $M\setminus\Omega$ is also an optimal domain, $M\setminus\Omega$ is also connected.  If $M=S^2$ then the Jordan curve theorem implies that $\Omega$ is simply connected.
\end{proof}

Another ingredient in our proof will be the following statement concerning the behaviour of the isoperimetric profile near $\xi=0$:

\begin{theorem}\label{thm:bdrybehaviour}
Let $M$ be a smooth compact Riemannian surface.  Then the isoperimetric profile satisfies 
$$
h(\xi)= \sqrt{4\pi|M|\xi}-\frac{|M|^{3/2}\sup_M K}{4\sqrt{\pi}}\xi^{3/2}+O(\xi^2)\quad\text{as}\quad\xi\to 0.
$$
\end{theorem}

Note that since $h(\xi)=h(1-\xi)$ this controls $h$ near both endpoints.

\begin{proof}
Small geodesic balls about any point $p$ are admissible in the definition of $h$.  The result of \cite{Gray}*{Theorem 3.1} gives $|B_r(p)|=\pi r^2\left(1-\frac{K(p)}{12}r^2+O(r^4)\right)$ and $|\partial B_r(p)|=2\pi r\left(1-\frac{K(p)}{6}r^2+O(r^4)\right)$.  The upper bound follows since $|\partial B_r(p)|\geq h(|B_r(p)|/|M|)$.

To prove the lower bound, first choose $\xi$ sufficiently small to ensure that $h(\xi)$ is much smaller than the injectivity radius of $M$.  Then the optimal region $\Omega$ corresponding to $\xi$ lies inside a geodesic ball about some point $p$, and the metric for $g$ in exponential coordinates is equal to the metric for the sphere of constant curvature $K(p)$ times a factor $(1+O(|\partial\Omega|^3))$.  The isoperimetric inequality for the sphere gives
\begin{align*}
h\left(\xi\right)&=|\partial\Omega|\\
&\geq \sqrt{4\pi|\Omega|}\left(1-\frac{K(p)|\Omega|}{8\pi}+O(|\partial\Omega|^3)\right)\\
&=\sqrt{4\pi|M|\xi}\left(1-\frac{K(p)|M|\xi}{8\pi}+O(\xi^{3/2})\right)\\
&\geq \sqrt{4\pi|M|\xi} - \frac{|M|^{3/2}\sup_MK}{4\sqrt{\pi}}\xi^{3/2}+O(\xi^2).
\end{align*}
\end{proof}

\section{A comparison theorem for the isoperimetric profile}

We begin with a technical statement comparing the isoperimetric profile of a solution of normalized Ricci flow with a subsolution of a certain parabolic equation.  Below we will see how to construct such subsolutions from particular solutions of Ricci flow.

\begin{theorem}\label{thm:comparison}
Let $\varphi:\ (0,1)\times[0,\infty)\to\RR$ be a smooth function which is strictly concave and positive for each $t\geq 0$, has $\limsup_{x\to 0}\frac{\varphi(\xi,t)}{4\pi\sqrt{\xi}}<1$ for every $t$, and satisfies
$$
\frac{\partial\varphi}{\partial t} < \frac{\varphi^2\varphi''-\varphi(\varphi')^2}{(4\pi)^2} + \varphi +\varphi'(1-2\xi)
$$
on $(0,1)\times[0,\infty)$.  Then for any smooth solution $(M,g(t))$ of the normalized Ricci flow on the two-sphere which satisfies $h_{g(0)}(\xi) > \varphi(\xi,0)$ for every $\xi\in(0,1)$, we have
$h_{g(t)}(\xi)\geq \varphi(\xi,t)$ for every $(\xi,t)\in(0,1)\times[0,\infty)$.
\end{theorem}

\begin{proof}
 We note that $h_{g(t)}(\xi)$ varies continuously in $t$ and in $\xi$, and that for $\xi$ sufficiently close to either $0$ or $1$ we have $h_{g(t)}(\xi)>\varphi(\xi,t)$ by Theorem \ref{thm:bdrybehaviour}.  Therefore if the claimed inequality does not hold for all time, there exists $t_0>0$ such that $h_{g(t)}(\xi,t)>\varphi(\xi,t)$ for all $\xi\in(0,1)$ and all $t\in[0,t_0)$, while $h_{g(t_0)}(\xi)\geq \varphi(\xi,t_0)$ for all $\xi\in(0,1)$ and there exists $\xi_0\in(0,1)$ with $h_{g(t_0)}(\xi_0)=\varphi(\xi,t_0)$.
 
Let $\Omega_0\subset M$ be the region in $M$ of area $4\pi\xi_0$ for which $|\partial\Omega|=\varphi(\xi_0,t_0)$, and let $\gamma_0=\partial\Omega_0$.  Since $h_{g(t)}(\xi)\geq \varphi(\xi,t)$ for $t\leq t_0$ and all $\xi$, we have
$$
|\gamma_0|_{g(t)}\geq \varphi\left(\frac{|\Omega_0|_{g(t)}}{4\pi},t\right)
$$
for $0\leq t\leq t_0$, and equality holds when $t=t_0$.  Since both sides of this equation are differentiable in $t$, it follows that under normalized Ricci flow,
\begin{equation}\label{eq:timeineq1}
\frac{\partial}{\partial t}|\gamma_0|_{g(t)}\Big|_{t=t_0}\leq \frac{\partial\varphi}{\partial t}+\frac{1}{4\pi}\varphi'(\xi_0,t_0)\frac{\partial}{\partial t}|\Omega_0|_{g(t)}\Big|_{t=t_0}.
\end{equation}
The time derivative on the left can be computed as follows:
$$
\frac{\partial}{\partial t}|\gamma_0|_{g(t)} = \frac{\partial}{\partial t}\int_{\gamma_0}\sqrt{g_t(\gamma_u,\gamma_u)}\,du = -\int_{\gamma_0}(K-1)ds=-\int_{\gamma_0}K\,ds + \varphi,
$$
where $ds$ is arc-length along $\gamma_0$.
Since $|\Omega_0|_{g(t)}=\int_{\Omega_0}\sqrt{\frac{\det g(t)}{\det g(t_0)}}\,d\mu(g(t_0))$, we find
$$
\frac{\partial}{\partial t}|\Omega_0|_{g(t)}\Big|_{t=t_0} = -2\int_{\Omega_0}(K-1)d\mu(g(t_0))
$$
Using the result of Theorem \ref{thm:curvetopology} we can apply the Gauss-Bonnet theorem, yielding
$$
\frac{\partial}{\partial t}|\Omega_0|_{g(t)}\Big|_{t=t_0} = 2|\Omega_0|-2\left(2\pi-\int_{\gamma_0}k\,ds\right)=8\pi\xi_0-4\pi+2\int_{\gamma_0}k\,ds,
$$
were $k$ is the geodesic curvature of the curve $\gamma_0$.  Thus the inequality \eqref{eq:timeineq1} becomes
\begin{equation}\label{eq:timeineq}
-\int_{\gamma_0}K\,ds + \varphi\leq \frac{\partial\varphi}{\partial t}+\frac{\varphi'}{4\pi}\left(8\pi\xi_0-4\pi+2\int_{\gamma_0}k\,ds\right).
\end{equation}

Now we observe that for any smooth family of domains $|\Omega(s)|$ in $M$ with $\Omega(0)=\Omega_0$, we have (at $t=t_0$) $|\partial\Omega(s)|\geq h\left(\frac{|\Omega(s)|}{4\pi}\right)\geq\varphi
\left(\frac{|\Omega(s)|}{4\pi},t_0\right)$, with equality at $s=0$.  It follows that the derivative with respect to $s$ of the difference is zero at $s=0$, and the second derivative is non-negative.  We consider for an arbitrary smooth function $\eta$ along $\gamma_0$ the smooth family of domains with boundary defined by $\gamma_s=\{\exp_z(s\eta(z)\vec{n}(z):\ z\in\gamma_0\}$ where $\vec{n}(z)$ is the unit normal vector to $\gamma_0$ at $z$ which points out of $\Omega_0$.  Then $\gamma_s$ is a smooth embedded closed curve for sufficiently small $s$, and we compute
\begin{equation*}
\frac{d}{ds}|\Omega(s)|\Big|_{s=0} = \int_{\gamma_0}\eta\,ds,
\end{equation*}
while 
\begin{equation*}
\frac{d}{ds}|\gamma_s|\Big|_{s=0} = \int_{\gamma_0}k\eta\,ds.
\end{equation*}
Therefore we have
$$
0=\int_{\gamma_0} k\eta\,ds - \frac{\varphi'}{4\pi}\int_{\gamma_0} \eta\,ds.
$$
Since $\eta$ is arbitrary, this implies that $k=\frac{\varphi'}{4\pi}$ along $\gamma_0$.  Now we consider the special case where $\eta\equiv 1$, in which case the motion of the boundary is at unit speed and we compute (noting $|\gamma_0|=\varphi$)
$$
\frac{d^2}{ds^2}|\Omega(s)|\Big|_{s=0} = \int_{\gamma_0}k\,ds = \frac{\varphi'}{4\pi}|\gamma_0|=\frac{\varphi\varphi'}{4\pi},
$$
and since $\frac{d}{ds}|\gamma_s| = \int_{\gamma_s}k\,ds = 2\pi-\int_{\Omega(s)}K\, d\mu$ for each $s$, we also have
$$
\frac{d^2}{ds^2}|\gamma_s|\Big|_{s=0} = -\int_{\gamma_0}K\,ds.
$$
Thus we have
\begin{equation}\label{eq:spaceineq}
0\leq -\int_{\gamma_0}K\,ds-\frac{\varphi(\varphi')^2}{(4\pi)^2}-\frac{\varphi^2\varphi''}{(4\pi)^2}.
\end{equation}

 Now we put the inequalities \eqref{eq:timeineq} and \eqref{eq:spaceineq} together, yielding
 \begin{align}\label{eq:ParabIneq}
 \frac{\partial\varphi}{\partial t}&\geq -\int_{\gamma_0}K\,ds+\varphi+\varphi'(1-2\xi_0)-2\left(\frac{\varphi'}{4\pi}\right)^2\varphi\notag\\
 &\geq \frac{\varphi(\varphi')^2}{(4\pi)^2}+\frac{\varphi^2\varphi''}{(4\pi)^2}+\varphi+\varphi'(1-2\xi_0)-2\frac{\varphi(\varphi')^2}{(4\pi)^2}\notag\\
 &=\frac{\varphi^2\varphi''-\varphi(\varphi')^2}{(4\pi)^2} + \varphi +\varphi'(1-2\xi_0).
 \end{align}
 This contradicts the strict inequality assumed in the theorem, and completes the proof.
 \end{proof}
 
 Now we proceed to the construction of solutions of the required differential inequality.  Our first step is to construct from any axially symmetric solution of normalized Ricci flow on the two-sphere a solution of the corresponding differential equation.  The construction is built from metrics of the form
$$
\tilde g = e^{2u(\phi)}\bar g
$$
where $\bar g$ is the standard metric on $S^2$, $u$ is a smooth even $\pi/2$-periodic function, and $\phi\in[0,\pi]$ is the angle defined by
 $$
 \cos\phi(x,y,z) = z.
 $$
 Given such a metric, the area of the spherical cap of angle $\phi$ is
 $$
 A_u(\phi) = 2\pi\int_0^\phi e^{2u(s)}\sin s\,ds
 $$
 and the length of its perimeter is
 $$
 L_u(\phi) = 2\pi e^{u(\phi)}\sin\phi.
 $$
 We require that the total area $A_u(\pi)$ of $\tilde g$ is $4\pi$.
 Observing that $A_u$ is a strictly increasing function, we have that $A_u$ has a well-defined
inverse from $[0,4\pi]$ to $[0,\pi]$ which is smooth on the interior.  We define
$$
\varphi_u(\xi):= L_u\circ A_u^{-1}(4\pi\xi).
$$

Under the normalized Ricci flow, the function $u$ evolves according to the equation
\begin{equation}\label{eq:2DRFu}
\frac{\partial u}{\partial t} = e^{-2u}\left(\frac{\partial^2u}{\partial\phi^2}+\cot\phi\frac{\partial u}{\partial\phi}\right)+1-e^{-2u}.
\end{equation}
 
We prove the following:

\begin{theorem}\label{thm:modelcase}
If $u$ is a solution of Equation \eqref{eq:2DRFu}, then the function $\varphi(\xi,t):=\varphi_{u(t)}(\xi)$
satisfies
\begin{equation}\label{eq:equality}
\frac{\partial\varphi}{\partial t} = \frac{\varphi^2\varphi''-\varphi(\varphi')^2}{(4\pi)^2} + \varphi +\varphi'(1-2\xi).
\end{equation}
Furthermore, $\varphi(\xi,t) = 4\pi\sqrt{\xi}\left(1-\frac12\left(1-2u''(0)\right)\xi+O(\xi^2)\right)$ as $\xi\to 0$.
\end{theorem}

\begin{proof}
By definition, the function $\varphi_{u(t)}(\xi)$ gives the length of the boundary of the spherical cap about the north pole which has area $4\pi\xi$ in the metric $g(t)=e^{2u(t)}\bar g$.  Therefore by construction for the family
of spherical caps $\Omega(s)=\{\phi\leq s\}$ we have the identity $|\partial\Omega(s)|_{g(t)}=\varphi
\left(\frac{|\Omega(s)|_{g(t)}}{4\pi},t\right)$ for every $s$ and every $t$.  Hence the time derivative of this difference is zero, and the second derivative with respect to unit speed normal motion is also zero, for each $s$ and $t$.  That is, equality holds in \eqref{eq:timeineq} and \eqref{eq:spaceineq}, and rearrangement gives the Theorem (note that by symmetry $k$ is constant on each of the curves $\partial\Omega_s$, so the vanishing of the first derivative with respect to unit speed motion implies $k=\frac{\varphi'}{4\pi}$ as required).  The asymptotic behaviour follows from the definitions of $L_u$ and $A_u$.
\end{proof}

We remark that the function $\varphi$ will not in general be the isoperimetric profile of the solution of normalized Ricci flow from which it is constructed.  For completeness we mention the following 
criterion for when this is the case, although it is not needed for our main application:

\begin{proposition}[\cite{Ritore}*{Theorem 3.5}]
If the curvature $\tilde K$ of $\tilde g_0=e^{2u_0}\bar g$ is positive and is decreasing in $\phi$ on the interval $(0,\pi/2)$, then $\varphi_{u(t)}(\xi)$ is the isoperimetric profile of $\tilde g(t)$ for every $t\geq 0$.
\end{proposition}

Finally, in order to apply Theorem \ref{thm:comparison}, we need to modify the functions constructed above to obtain strict inequalities required both in the differential inequality and at the boundary, as well as the strict concavity of $\varphi$.

The main requirement to accomplish this is that $\tilde g$ have positive curvature (note that if this is true initially, then it remains true for positive times by the maximum principle).  In this case the equality in \eqref{eq:spaceineq} gives
$$
0=-\int_{\gamma_0}\tilde K\,ds -\frac{\varphi(\varphi')^2+\varphi^2\varphi''}{(4\pi)^2},
$$
so that 
\begin{equation}\label{eq:concave}
\varphi'' = -(4\pi)^2 K - \frac{(\varphi')^2}{\varphi}<0.
\end{equation}

\begin{theorem}\label{thm:comparison2}
Let $\tilde g(t)$ be a rotationally symmetric solution of normalized Ricci flow on the two-sphere with positive Gauss curvature, and let $\varphi$ be the corresponding solution of equation \eqref{eq:equality}.
Let $g(t)$ be any smooth solution of normalized Ricci flow on the two-sphere with $h_{g(0)}(\xi)\geq\varphi(\xi,0)$ for all $\xi\in(0,1)$.  Then $h_{g(t)}(\xi)\geq \varphi(\xi,t)$ for all $\xi\in(0,1)$ and all $t$ in the interval of existence of $g$ and $\tilde g$.
\end{theorem}

\begin{proof}
For any $\varepsilon\in(0,1)$ define $\varphi_\varepsilon(\xi,t)=(1-\varepsilon)\varphi(\xi,t)$.  Then $\varphi(.,t)$ is strictly concave for each $t$, we have $h_{g(0)}(\xi)>\varphi_\varepsilon(\xi,0)$ for all $\xi$, $\limsup_{\xi\to 0}\frac{\varphi_\varepsilon(\xi,t)}{4\pi\sqrt{\xi}}=1-\varepsilon$ for each $t$, and 
\begin{align*}
\frac{\partial\varphi_\varepsilon}{\partial t} &- \frac{\varphi_\varepsilon^2\varphi_\varepsilon''-\varphi_\varepsilon(\varphi_\varepsilon')^2}{(4\pi)^2} - \varphi_\varepsilon -\varphi_\varepsilon'(1-2\xi)\\
&=\frac{\varepsilon(1-\varepsilon)(2-\varepsilon)}{(4\pi)^2}
\left(\varphi^2\varphi''-\varphi(\varphi')^2\right)\\
&<0
\end{align*}
by the estimate \eqref{eq:concave}.  Thus Theorem \ref{thm:comparison} applies to yield
$h_{g(t)}(\xi)\geq \varphi_{\varepsilon}(\xi)$
for all $\xi$ and $t$, and the result follows after letting $\varepsilon$ approach zero.
\end{proof}
 
\section{Comparison with the Rosenau solution}

The Rosenau solution is an explicit axially symmetric solution of (normalized) Ricci flow on the two-sphere.  The metric is given by $g(t) = u(x,t)(dx^2+dy^2)$, where $(x,y)\in\RR\times[0,4\pi]$, where
$$
u(x,t) = \frac{\sinh(e^{-2t})}{2e^{-2t}\left(\cosh(x)+\cosh(e^{-2t})\right)}.
$$
This extends to a smooth metric on the two-sphere at each time with area $4\pi$, which evolves according to the normalized Ricci flow equation \eqref{eq:nRF}.  A direct computation gives the
corresponding function constructed in Theorem \ref{thm:modelcase} to be
\begin{equation}\label{eq:Rosenauprofile}
\varphi(\xi,t) = 4\pi\sqrt{\frac{\sinh(\xi e^{-2t})\sinh((1-\xi)e^{-2t})}{\sinh(e^{-2t})e^{-2t}}}.
\end{equation}
In particular the asymptotic behaviour  is given by 
$$
\varphi(\xi,t)=4\pi\sqrt{\xi}\left(1-\frac12e^{-2t}\coth(e^{-2t})\xi+O(\xi^2)\right)\qquad\text{ as\ }\xi\to 0.
$$
In order to apply theorem \ref{thm:comparison2} we need to be able to establish the required inequality at the initial time:

\begin{theorem}\label{thm:initialrosenaucomp}
For any smooth metric $g_0$ on the two-sphere, there exists $t_0\in\RR$ such that
$h_{g_0}\geq \varphi(.,t_0)$, where $\varphi$ is as given in Equation \eqref{eq:Rosenauprofile}.
\end{theorem}

\begin{proof}
$\varphi$ is continuous on $[0,1]\times\RR$, and strictly increasing in $t$ for each $\xi\in(0,1)$.  Also, we have $\lim_{t\to\infty}\varphi(\xi,t) = 4\pi\sqrt{\xi(1-\xi)}$ and $\lim_{t\to-\infty}\varphi(\xi,t)=0$.  Therefore at values of $\xi\in(0,1)$ where $h_{g_0}(\xi)<4\pi\sqrt{\xi(1-\xi)}$ there is a unique $t(\xi)\in\RR$ such that
$h_{g_0}(\xi)=\varphi(\xi,t(\xi))$.  Furthermore, the asymptotic behaviour of $\varphi$ as $\xi\to 0$ and 
Theorem \ref{thm:bdrybehaviour} imply that $e^{-2t_*}\coth(e^{-2t_*})=\sup_MK$, where $t_*=\lim_{\xi\to 0}t(\xi)$.  The left-hand side has decreases from infinity to $1$ as $t_*$ increases from $-\infty$ to $\infty$, and so there is a unique solution for $t_*$.  That is, $t(\xi)$ extends to a continuous function from $[0,1]$ to $\RR\cup\{\infty\}$, and so is bounded below.  Choosing $t_0=\inf\{t(\xi):\ \xi\in(0,1)\}$ gives the result.
\end{proof}

This, together with Theorems \ref{thm:comparison2} and \ref{thm:bdrybehaviour} yields the following:

\begin{theorem}\label{thm:Curvest}
For any smooth solution of normalized Ricci flow on the two-sphere, there exists $t_0\in\RR$ such that
$$
K(x,t)\leq \coth\left(\E^{-2(t+t_0)}\right)\E^{-2(t+t_0)}
$$
for all points $x$ and all $t$ in the interval of existence of the solution.  It follows that $K(x,t)\leq \sqrt{1+\E^{-4(t+t_0)}}\leq 1+\frac12 \E^{-4(t+t_0)}$.
\end{theorem}

\section{Convergence to constant curvature}

Since the curvature evolves by $\partial_tK = \Delta K + K(K-1)$,  we have the lower bound $K\geq -\frac{1}{\E^t-1}$ for $t>0$.  Thus by Theorem \ref{thm:Curvest} the curvature remains bounded and the solution exists for all positive times, and we have bounds of the form
\begin{equation}\label{eq:highreg}
|\nabla^{(k)}K|^2\leq C(k,t_0)(1+t^{-k})
\end{equation}
 for every $k$, by a maximum principle boot-strapping argument as described in \cite{HamFSRF}*{Section 7}.
We show that the curvature converges in $L^1$ to $1$:  By the Gauss-Bonnet theorem and Theorem \ref{thm:Curvest} we have
$$
0=\int_MK-1\,d\mu \leq -\int_{K\leq 1}|K-1|\,d\mu + 2\pi\E^{-4(t+t_0)}.
$$
Rearranging gives $\int_{K\leq 1}|K-1|\,d\mu\leq 2\pi\E^{-4(t+t_0)}$, and so
\begin{equation}\label{eq:curvdecay}
\int_M|K-1|\,d\mu\leq 4\pi\E^{-4(t+t_0)}.
\end{equation}  

The control on the isoperimetric profile bounds the isoperimetric constant $\frak{I} = \inf\left\{
\frac{|\partial\Omega|^2}{\min\{|\Omega|,|M\setminus\Omega|\}}:\ \Omega\subset M\right\}$, and so also controls the constant in the Sobolev inequality (see \cite{Chavel}*{Theorem 12, p.111}).  This also bounds the constants in the Gagliardo-Nirenberg inequalities (see the proof in \cite{Aubin}*{page 93}), so
$$
\|\nabla^kK\|_{\infty}\leq C(k,m,t_0)\|K-1\|_{L^1}^{\frac{m-k}{m+2}}\|\nabla^mK\|_{\infty}^{\frac{k+2}{m+2}}\leq C(k,\varepsilon,t_0)\E^{-(4-\varepsilon)t},
$$
for any $\varepsilon>0$ for $t\geq 1$, using the estimates \eqref{eq:curvdecay} and \eqref{eq:highreg} with $m$ large enough for given $\varepsilon>0$.
The convergence of the metric to a constant curvature metric now follows:  The metric remain comparable to the initial metric and converges uniformly to a limit metric, since for any nonzero $v\in TM$, 
$$
\left|\frac{\partial}{\partial t}\log g(v,v)\right|=2|K-1|\leq C\E^{-2t}
$$
which is integrable.  Convergence in $C^\infty$ follows as in \cite{Ham82}*{Section 17}.

\end{document}